\documentclass[12pt]{amsart}
\usepackage[top=1in, bottom=1in, left=1in, right=1in]{geometry}
\usepackage{dsfont}            
\usepackage{amssymb}  
\usepackage[utf8]{inputenc}
\usepackage{mathtools}
\usepackage{amsthm}  
\usepackage{indentfirst}  
\usepackage{amsmath}
\usepackage{graphicx}
\usepackage{times}
\usepackage{color}
\usepackage[shortlabels]{enumitem}
\usepackage[colorlinks=true]{hyperref}
\hypersetup{colorlinks=true, citecolor=green, linkcolor=green, urlcolor=cyan, filecolor=magenta}

\numberwithin{equation}{section}

\newtheorem{theorem}{Theorem}[section]
\newtheorem{proposition}[theorem]{Proposition}  
\newtheorem{corollary}[theorem]{Corollary}  
\newtheorem{lemma}[theorem]{Lemma}

\newtheoremstyle{claim} 
{1em}                        
{1em}                        
{}                           
{}                           
{\bfseries}                  
{.}                          
{.5em}                       
{}  

\theoremstyle{claim}    
\newtheorem{claim}{Claim}

\theoremstyle{remark}
\newtheorem{remark}{Remark}[section]
\newtheorem*{remark*}{Remark}

\theoremstyle{definition}

\newcommand{\N}{\mathbb{N}}
\newcommand{\Z}{\mathbb{Z}}
\newcommand{\Q}{\mathbb{Q}}
\newcommand{\R}{\mathbb{R}}
\newcommand{\C}{\mathbb{C}}

\newcommand{\CF}{\mathcal{F}}

\newcommand{\CP}{\mathcal{P}}

\DeclareMathOperator{\supp}{supp}

\newcommand{\dee}{\mathrm{d}}

\newcommand{\abs}[1]{\left\lvert #1 \right\rvert}

\newcommand{\bs}\boldsymbol{}
\renewcommand{\geq}{\geqslant}
\renewcommand{\leq}{\leqslant}

\renewcommand{\tilde}{\widetilde}

\renewcommand{\Re}{{\rm Re}}

\definecolor{blue}{rgb}{.2,.6,.75}
\definecolor{green}{rgb}{.4,.7,.4}

\begin{document}

\title[On Multiplicative Functions with Small Partial Sums]{On Multiplicative Functions with Small Partial Sums}

\author{Stelios Sachpazis}
\address{D\'epartement de math\'ematiques et de statistique\\
Universit\'e de Montr\'eal\\
CP 6128 succ. Centre-Ville\\
Montr\'eal, QC H3C 3J7\\
Canada}
\email{stelios.sachpazis@umontreal.ca}


\date{\today}

\begin{abstract}
In analytic number theory, several results make use of information regarding the prime values of a multiplicative function in order to extract information about its averages. Examples of such results include Wirsing's theorem and the Landau\,-\,Selberg\,-\,Delange method. In this paper, we are interested in the opposite direction. In particular, we prove that when $f$ is a suitable divisor-bounded multiplicative function with small partial sums, then $f(p)\approx-p^{i\gamma_1}-\ldots-p^{i\gamma_m}$ on average, where the $\gamma_j$'s are the imaginary parts of the zeros of the Dirichet series of $f$ on the line $\Re(s)=1$. This extends a result of Koukoulopoulos and Soundararajan and it builds upon ideas coming from previous work of Koukoulopoulos for the case where $|f|\leq 1$.
\end{abstract}

\maketitle


\section{Introduction}
 
Let $D$ be a fixed positive integer and let $f$ be a multiplicative function such that $|f|\leq\tau_D$. Recent work of Granville and Koukoulopoulos \cite{blsd} implies that if $f(p)\approx v$ on average, in the sense that

\begin{equation}\label{clsd}
\sum_{p\leq x}f(p)\log p=vx+O_M\left(\frac{x}{(\log x)^M}\right)
\end{equation}

\vspace{2mm}

\noindent
for all $x\geq 2$ and $M>0,$ then

\begin{equation}\label{lsd}
\sum_{n\leq x}f(n)=\frac{c(f,v)}{\Gamma(v)}x(\log x)^{v-1}+O_{D,M,v}\left(x(\log x)^{\Re(v)-2}\right)
\end{equation}

\vspace{2mm}

\noindent
for some constant $c(f,v).$ This result builds on earlier work of Landau-Selberg-Delange \cite{del,lan1,lan2,sel} and of Wirsing \cite{wir1,wir2}.

Since $|f|\leq \tau_D$, condition (\ref{clsd}) implies that $|v|\leq D$. So, if $v\notin \Z_{\leq 0}$ and $c(f,v)\neq 0$, then the ``main term'' of the asymptotic formula (\ref{lsd}) does not vanish and it thus determines the asymptotic behaviour of the partial sums of $f$. More precisely, it is then true that

\begin{equation*}
\left|\sum_{n\leq x}f(n)\right|\sim \abs{\frac{c(f,v)}{\Gamma(v)}}\frac{x}{(\log x)^{1-\Re({v})}} \quad \text{as} \quad x\rightarrow \infty.
\end{equation*}

\vspace{3mm}

\noindent
Of course, $\Re(v)\geq -|v|\geq -D$, and this means that, asymptotically, the size of the partial sums $\sum_{n\leq x}f(x)$ can get no smaller than $x/(\log x)^{D+1}.$ 

Hence, if $c(f,v)\neq 0$, then by further assuming that

\begin{equation}\label{cond}
\sum_{n\leq x}f(n)\ll \frac{x}{(\log x)^A}\quad \text{for some fixed}\quad \!\!A>D+1,
\end{equation}

\vspace{2mm}

\noindent
we force the complex number $v$ to be a non-positive integer of $\{0,-1,\ldots,-D\}.$

The trivial upper bound for $\sum_{n\leq x}f(n)$ is $x(\log x)^{D-1}$, the upper bound that one can get by replacing $f$ by $\tau_D.$ Therefore, the above reasoning claims that when $c(f,v)\neq 0$ and $f$ has small averages, that is (\ref{cond}) is satisfied, there exists only a finite number of possibilities for the average size of $f(p)$. However, the respective argument was developed under condition (\ref{clsd}), a condition that guarantees a priori knowledge on the average size of $f$ on the primes. Consequently, we might start wondering whether it is possible to draw specific conclusions about the average size of $f(p)$ when (\ref{cond}) holds without assuming (\ref{clsd}). The goal of the present paper is to provide an answer to this question. To this end, in the last section we prove the following theorem, which is our main result.

\vspace{2mm}

\begin{theorem}\label{mthm}
Fix a positive integer $D$ and a real number $A>D+1$. Let also $f$ be a multiplicative function such that $|\Lambda_f|\leq D\cdot \Lambda$, where $\Lambda_f$ is the unique arithmetic function defined through the relation $f\cdot \log=\Lambda_f\ast f$ and $\Lambda$ is the von Mangoldt function. Assume further that

$$\sum_{n\leq w}f(n)\ll \frac{w}{(\log w)^A}\quad \text{for all}\quad \!w\geq 2.$$

\vspace{1mm}

\noindent
Then there exists a multiset $\Gamma$ of $m$ real numbers with $m\leq D$ and such that 

\begin{equation}\label{mainin}
\Big|\sum_{p\leq x}\Big(f(p)+\sum_{\gamma\in \Gamma}p^{i\gamma}\Big)\log p\Big|\leq O_{f,T}\left(\frac{x(\log\log x)^{D+m}}{(\log x)^{\min\{1,A-D-1\}/2}}\right)+O_{\Gamma}\!\left(\frac{x(\log T)^{D+m}}{\sqrt{T}}\right)\!\!,
\end{equation}

\vspace{2mm}

\noindent
for all $x\geq 3$ and any $T\geq 1$. The implied constants depend also on $D$ and $A$.
\end{theorem}

\vspace{1mm}

\begin{corollary}
Under the assumptions of Theorem \ref{mthm}, there exists a multiset $\Gamma$ of at most $D$ real numbers such that 

\begin{equation}\label{limit}
\lim_{x\to+\infty}\frac{1}{x}\sum_{p\leq x}\Big(f(p)+\sum_{\gamma\in \Gamma}p^{i\gamma}\Big)\log p=0.
\end{equation}
\end{corollary}

\vspace{2mm}

So, Theorem \ref{mthm} predicts that $f(p)\approx-\sum_{\gamma\in\Gamma}p^{i\gamma}$ on average. It turns out that $\Gamma$ is the multiset of the ordinates of the zeros of the Dirichlet series of $f$ on the line $\Re(s)=1$.

Results like Wirsing's theorem \cite{wir1,wir2} and the Landau-Selberg-Delange method \cite{del,lan1,lan2,sel} require information about the prime values of a multiplicative function and they give us back information about its averages. Theorem \ref{mthm} does the converse. It uses the fact that the averages of an appropriate multiplicative function are small and it returns information about its structure on primes.

In 2013, Koukoulopoulos \cite{oldk} covered the case where $D=1$ and $A>2$. Seven years later, in 2020, in joint work with Soundararajan \cite{kousou}, they proved an analogue of Theorem \ref{mthm} when $D$ is any fixed positive integer and $A>D+2$. In the special case where the Dirichlet series of $f$ has a single root of multiplicity $D$ on the line $\sigma=1$, they proved the same result in the range $A>D+1$. Theorem \ref{mthm} extends both of their results, since it takes the full range $A>D+1$ into account and makes no assumptions on the zeros of the Dirichlet series of $f$.

\begin{remark}
If the hypothesis $A>D+1$ is replaced by the inequality $A<D+1$, then the prediction of Theorem \ref{mthm} is false and we can verify this with a counterexample. Let $\kappa\in\R\setminus\Q$ be such that $A-1<\kappa<D$. Then $|\Lambda_{\tau_{-\kappa}}|=\kappa\cdot\Lambda\leq D\cdot\Lambda$ and $\sum_{n\leq x}\tau_{-\kappa}(n)\ll_{\kappa}x/(\log x)^{\kappa+1}\leq x/(\log x)^A$ for all $x\geq 3$, meaning that $\tau_{-\kappa}$ would meet the conditions of Theorem \ref{mthm} with $A<D+1$. The Dirichlet series of $\tau_{-\kappa}$ is $\zeta^{-\kappa}(s)$ and it has no zeros on $\Re(s)=1$, because $\kappa$ is irrational and this implies that $\zeta^{-\kappa}$ has a singularity at $s=1$ which is not a pole. So, the multiset $\Gamma$ is empty in the case of $\tau_{-\kappa}$. Consequently, under the assumption $A<D+1$, for the function $\tau_{-\kappa}$, Theorem \ref{mthm} would yield that $-\kappa=\tau_{-\kappa}(p)\approx 0$ on average, which is absurd.	
\end{remark}

\begin{remark}
Note that in Theorem \ref{mthm} one can restrict the sum $\sum_{\gamma\in\Gamma}p^{i\gamma}$ to the ordinates $|\gamma|\leq T$. Indeed, for $|\gamma|>T$, we have that

$$\sum_{p\leq x}p^{i\gamma}\log p=\frac{x^{1+i\gamma}}{1+i\gamma}+O\Big(\frac{x}{\log x}\Big)\ll\frac{x}{T}+\frac{x}{\log x}.$$
\end{remark}

\vspace{1mm}

\subsection*{Notation and definitions}
For notational ease, throughout the rest of the text we supress any dependence of the implied constants on $D$ and $A$. We are doing this because there are many such dependencies. For an integer $n>1$, we denote the smallest prime factor of $n$ by $P^-(n)$. For $n=1$, we define $P^-(1)=+\infty$. Moreover, in this paper, the lowercase greek letter $\phi$ denotes the Euler totient function and the arithmetic functions $\omega$ and $\Omega$ are defined, as usual, through the relations $\omega(n)=\sum_{p\mid n}1$ and $\Omega(n)=\sum_{p^a\| n}a$. An arithmetic function which is bounded by a divisor function $\tau_{k}$ will be called divisor-bounded. Finally, we will be denoting the Dirichlet series of an arithmetic function $f$ by

$$L(s,f)=\sum_{n\geq 1}\frac{f(n)}{n^s},$$

\vspace{1mm}

\noindent
provided that the series of the right-hand side converges.

\vspace{1mm}

\subsection*{Acknowledgements}
The author would like to thank his advisor, Dimitris Koukoulopoulos, for his constant guidance and support during the making of this work. Their many fruitful mathematical discussions on the subject proved to be very important for the development of this paper. Last, but not least, he would also like to wholeheartedly express his gratitude towards the Stavros Niarchos Scholarships Foundation for all the financial support that provides to him for his doctoral studies.

\vspace{2mm}

\section{The Ideas Behind the Proof of Theorem \ref{mthm}}

In this section we sketch and motivate the proof of Theorem \ref{mthm}. But, before doing so, we give some necessary notation. 

First, as in the statement of Theorem \ref{mthm}, for a multiplicative function $g$, we will be using the notation $\Lambda_g$ for the unique arithmetic function defined through the relation

$$g\cdot \log=\Lambda_g\ast g.$$

\vspace{2mm}

\noindent
The function $\Lambda_g$ is supported on prime powers and in particular, $\Lambda_g(p)=g(p)\log p$ on the primes. It is also true that ${\Lambda}_{g\ast h}=\Lambda_g+\Lambda_h$ for any two multiplicative functions $g$ and $h$.

We now introduce two classes of multiplicative functions and make a few comments about the first one. Given an integer $D\in \mathbb{N}$ and a real number $A>0$, we define the sets

$$\begin{gathered}
\CF(D):=\{f:\N\rightarrow \C,\,f\,\text{multiplicative},\, |\Lambda_f|\leq D\cdot \Lambda\},\\
\mathcal{F}(D,A):=\left\{f\in \mathcal{F}(D),\ \sum_{n\leq x}f(n)\ll \frac{x}{(\log x)^A}\ \,\text{for all}\ x\geq 2\right\}\!.
\end{gathered}$$

\vspace{2mm}

Note that for a function $f\in\CF(D)$, we have that

$$-\frac{L'}{L}(s,f)=\sum_{n\geq 1}\frac{\Lambda_f(n)}{n^s}\quad \text{for}\quad \Re(s)>1$$

\vspace{2mm}

\noindent
and that the series convergences absolutely when $\Re(s)>1$.

The class $\CF(D)$ includes many important number-theoretic functions, like the M\"{o}bius function $\mu$, the Liouville function $\lambda$ and the generalized divisor functions $\tau_k$ for $k\leq D$. There are technical reasons that make the class $\CF(D)$ very convenient to work with. For example, if $f\in\CF(D),$ then $f^{-1}\in\CF(D)$, where $f^{-1}$ is the Dirichlet inverse of $f$, namely, the inverse of $f$ with respect to the Dirichlet convolution. Furthermore, if $f\in\CF(D)$, then both $f$ and $f^{-1}$ satisfy the inequalities 

$$|f|\leq {\tau}_D, \quad |f^{-1}|\leq {\tau}_D.$$

\vspace{3mm}

\noindent
These two last results are proved in \cite[Lemma 2.2]{kousou}.

Having now defined the classes $\CF(D)$ and $\CF(D,A)$, we continue by describing the ideas behind the proof of Theorem \ref{mthm}. For a multiset $\Gamma=\{\gamma_1,\ldots,\gamma_m\}$, we let

$$\tau_{\Gamma}(n)=\sum_{d_1\cdots d_m=n}d_1^{i\gamma_1}\cdots d_m^{i\gamma_m} \quad \text{for} \quad n\in\mathbb{N}.$$

\vspace{2mm}

\noindent
Then, for a function $f\in\CF(D,A)$, we consider the multiplicative function $f_{\Gamma}=f\ast\tau_{\Gamma}$. On the primes $p$, the values of the function $f_{\Gamma}$ are

$$f_{\Gamma}(p)=f(p) +\sum_{\gamma\in\Gamma}p^{i\gamma}.$$

\vspace{1mm}

\noindent
From now on, we focus on establishing a bound for the sums $\sum_{p\leq x}f_{\Gamma}(p)\log p$ for some suitable multiset $\Gamma$. Note that

\begin{equation}\label{id}
{\Lambda}_{f_{\Gamma}}(n)={\Lambda}_f(n)+\sum_{\gamma\in \Gamma}\Lambda(n)n^{i\gamma}.
\end{equation}

\vspace{1mm}

\noindent
Therefore, $|{\Lambda}_{f_{\Gamma}}|\leq (D+m)\cdot\Lambda$, as $f\in\CF(D)$. Moreover, since $\Lambda_{f_{\Gamma}}$ is supported on prime powers, and the contribution of the squares of primes, the cubes of primes and all the larger powers of primes in $\sum_{n\leq x}\Lambda_{f_{\Gamma}}(n)$ is at most $O_m(\sqrt{x})$, we have

$$\sum_{p\leq x}f_{\Gamma}(p)\log p=\sum_{n\leq x}\Lambda_{f_{\Gamma}}(n)+O_m(\sqrt{x}),$$

\vspace{2mm}

\noindent
and so our goal now is to estimate the sums of the right-hand side.

With $\Gamma$ being the multiset of the ordinates of the zeros of the continuous extension of $L(\cdot,f)$ on $\Re(s)\geq 1$, Koukoulopoulos and Soundararajan \cite{kousou} bounded the sums $\sum_{n\leq x}\Lambda_{f_{\Gamma}}(n)$ by applying a smoothed version of Perron's formula to the Dirichlet series of $-(L'/L)'(\cdot,f_{\Gamma})$. Their approach requires information on the size of $L(\cdot,f_{\Gamma})$ and of its two first derivatives. It turns out that this information is available once $A>D+2$. The wider range $A>D+1$ can only ensure bounds for $L(\cdot,f_{\Gamma})$ and $L'(\cdot,f_{\Gamma})$. However, if they had used a variant of Perron's formula for $(L'/L)(\cdot,f_{\Gamma})$, so that they avoided the presence of the second derivative $L''(\cdot,f_{\Gamma})$, the resulting bound at the end of their argument would have been off by one factor of $\log x$, making it trivial.

In the present paper, we circumvent the use of higher derivatives which are responsible for the condition $A>D+2$ in the work of Koukoulopoulos and Soundararajan \cite{kousou}. We do so 
by resorting to the recursiveness of the mean values of multiplicative functions. By this, we mean the identity

\begin{equation}\label{ref}
\sum_{n\leq x}f_{\Gamma}(n)\log n=\sum_{n\leq x}\Lambda_{f_{\Gamma}}(n)\sum_{d\leq x/n}f_{\Gamma}(d).
\end{equation}

\vspace{2mm}

\noindent
The idea that is about to be described is inspired by the work of Koukoulopoulos \cite{oldk}. By establishing an estimate for the sums $\sum_{n\leq x}f_{\Gamma}(n)$, we can apply partial summation to bound the sums of the left-hand side in (\ref{ref}). One might then make use of Dirichlet's hyperbola method to obtain an estimate for

$$\sum_{n\leq \sqrt{x}}f_{\Gamma}(n)\sum_{d\leq x/n}\Lambda_{f_{\Gamma}}(d).$$

\vspace{2mm}

\noindent
The summand corresponding to $n=1$ is $\sum_{d\leq x}\Lambda_f(d)$, which is the sum that we are trying to bound. The problem is that the next term equals $f(2)\sum_{d\leq x/2}\Lambda_f(d)$ and this term is expected to have roughly the same size as the ``main term'' $\sum_{d\leq x}\Lambda_f(d)$. In order for this obstacle to be overcome, sieve methods come into play. By combining sieves with the hyperbola method and the ``recursiveness'' of averages for the function $f_{\Gamma}\cdot \mathds{1}_{P^-(\cdot)>z}$, we aim to estimate

$$\sum_{\substack{n\leq \sqrt{x}\\P^-(n)>z}}\!\!\!f_{\Gamma}(n)\sum_{d\leq x/n}\Lambda_{f_{\Gamma}}(d).\,$$

\vspace{2mm}

\noindent
The summand for $n=1$ is again the sum $\sum_{d\leq x}\Lambda_f(d)$, but this time all the summands for $n\in(1,z]$ vanish. So, the problem that occurred before is now resolved. In addition, the summands with $n>z$ are supported on a set of density $\ll 1/\log z$.

With this approach we will bound $\sum_{d\leq x}\Lambda_{f_{\Gamma}}(d)$ in terms of its integrals

$$\int_{1}^x\frac{\left|\sum_{d\leq t}\Lambda_{f_{\Gamma}}(d)\right|}{t^2}\mathrm{d}t.$$

\vspace{2mm}

\noindent
We finally estimate these integrals by proceeding in the following way. After an application of the Cauchy-Schwarz inequality, we apply Parseval's theorem to bound the above integral by an integral involving the logarithmic derivative of $L(\cdot,f_\Gamma)$. We continue by splitting the integral in the ranges $|t|\leq T$ and $|t|>T$. In the first range we bound the integral trivially by using the continuity of $L(\cdot,f_{\Gamma})$ and $L'(\cdot,f_{\Gamma})$ and in the second one we apply a mean value theorem about $L^2$ norms of Dirichlet series.

\vspace{2mm}

\section{Auxiliary Results}

This section is a collection of lemmas, propositions and theorems that will be of use later in the text. We begin with a lemma \cite[Theorem 2.22]{tom} which is an immediate consequence of M\"{o}bius inversion.

\begin{lemma}\label{mobin}
Let $F:(0,+\infty)\rightarrow \C$ and $G:(0,+\infty)\rightarrow \C$ be two complex-valued functions such that $F(x)=G(x)=0$ for $x\in(0,1)$. Let $h$ be an arithmetic function which has an inverse under Dirichlet
convolution which we denote by $h^{-1}$. If 
	
$$G(x)=\sum_{n\leq x}h(n)F(x/n),$$
	
\noindent
then
	
$$F(x)=\sum_{n\leq x}h^{-1}(n)G(x/n).$$
\end{lemma}

\vspace{1mm}

\begin{proposition}\label{geo}
Let $\{a_n\}_{n\in\N}$ be a sequence of positive real numbers and let $k$ be a positive integer. For $y\geq 0$, we have that
	
$$\#\Big\{(\nu_1,\ldots,\nu_k)\in{\N}^k:\sum_{j=1}^ka_j\nu_j\leq y\Big\}\leq \frac{\big(y+\sum_{j=1}^ka_j\big)^k}{k!\prod_{j=1}^ka_j}.$$
\end{proposition}

\begin{proof}
See \cite[Theorem 3, p.~363]{tenen}.
\end{proof}

\vspace{1mm}

\begin{theorem}\label{fromb}
If $f$ is a multiplicative function such that $0\leq f\leq {\tau}_k$ for some positive integer $k$, then

$$\sum_{n\leq x}f(n)\ll_{k}x\exp\Big\{\sum_{p\leq x}\frac{f(p)-1}{p}\Big\}.$$
\end{theorem}

\begin{proof}
See \cite[Theorem 14.2]{dimb}.
\end{proof}

\vspace{1mm}

For the next lemma we introduce some notation which will appear again in Proposition \ref{corpr}: for $t\in\R$, set

$$V_t:=\exp\{100(\log (|t|+3))^{2/3}(\log \log (|t|+3))^{1/3}\}.$$

\vspace{1mm}

\begin{lemma}\label{dle}
Let $z\geq 2$ and $t\in \R$ such that $z\geq V_t$. For $x\geq z$, we have that 

$$\sum_{\substack{n\leq x\\P^-(n)>z}}\!\!\!n^{it}=\frac{x^{1+it}}{1+it}\prod_{p\leq z}\left(1-\frac{1}{p}\right)+O\left(\frac{x^{1-1/(30\log z)}}{\log z}\right)\!.$$
\end{lemma}

\begin{proof}
See \cite[Lemma 3.1]{oldk}.
\end{proof}

\vspace{0.1mm}

\begin{proposition}\label{dsp}
Let $f$ be a function of the class $\CF(D,A)$ with $A>D+1$.

\begin{enumerate}[(a)]

    \item The series $L^{(j)}(s,f)$ with $0\leq j<A-1$ 
    all converge uniformly in compact subsets of the region $\Re(s)\geq 1$.
	
    \item Counted with multiplicity, $L(s,f)$ has at most $D$ zeros on the line $\Re(s)=1$.
	
    \item Let $\Gamma$ denote the (possibly empty) multiset of ordinates $\gamma$ of zeros $1+i\gamma$ of $L(s,f)$. Let $\tilde{\Gamma}$ denote a (multi-)subset of $\Gamma$ and let $m_{\tilde{\Gamma}}$ denote the largest multiplicity of an element in $\tilde{\Gamma}$. The Dirichlet series
	
	$$L(s,f_{\tilde{\Gamma}})=L(s,f)\prod_{\gamma\in\tilde{\Gamma}}\zeta(s-i\gamma)$$
	
	\noindent
	and the series of derivatives $L^{(j)}(s,f_{\tilde{\Gamma}})$ for $1\leq j<A-m_{\tilde{\Gamma}}-1$ all converge uniformly in compact
	subsets of the region $\Re(s)\geq 1$.
\end{enumerate}
\end{proposition}

\begin{proof}
See Proposition 2.4 in \cite{kousou}.
\end{proof}

\vspace{1mm}

\begin{theorem}\label{1r}
Fix a positive integer $D$ and a real number $A>D+1$. If $f\in\CF(D,A)$ and $L(\cdot,f)$ has one zero $1+i\gamma$ of multiplicity $D$ on the line $\Re(s)=1$, then

\begin{equation*}
\sum_{p\leq x}|f(p)+Dp^{i\gamma}|\log p\ {\ll}_{f}\frac{x}{(\log x)^{(A-1-D)/2}}.
\end{equation*} 
\end{theorem}

\begin{proof}
See \cite[p.~12-13]{kousou}.
\end{proof}

\vspace{1mm}

The following result is due to Montgomery \cite[Theorem 3, p. 131]{mon} and it is a classic mean value theorem for Dirichlet series.

\vspace{1mm}

\begin{lemma}\label{mvt}
Let $A(s)=\sum_{n\geq 1}a_nn^{-s}$ and $B(s)=\sum_{n\geq 1}b_nn^{-s}$ be two Dirichlet series which converge for $\Re(s)>1.$ If $|a_n|\leq b_n$ for all $n\in \N$, then
	
$$\int_{-T}^T|A(\sigma+it)|^2\dee t\leq 3\int_{-T}^T|B(\sigma+it)|^2\dee t,$$

\vspace{2mm}
	
\noindent
for any $\sigma>1$ and any $T\geq 0$. 
\end{lemma}

\vspace{1mm}

The last result of this section is a standard sieve type lemma which may be found in \cite{dimb} as part (b) of Exercise 19.4.

\vspace{1mm}

\begin{lemma}\label{sieex}
Let $k,\kappa$ and $C$ be positive constants and let also $z\geq 2,\, \CP\subseteq \{p\leq z\}, P(y)=\prod_{p\in\CP, p\leq y}p$ for all $y\leq z$ and $u\geq u_{\kappa}:=1+2/(e^{0.53/\kappa}-1)$.  Suppose further that $\nu$ is a multiplicative function such that $0\leq \nu(p)<\min\{k,p\}$ for all $p\in\CP$ and

$$\prod_{\substack{p\in \CP,\\w_1<p\leq w_2}}\!\!\left(1-\frac{\nu(p)}{p}\right)^{-1}\!\!\leq \left(\frac{\log w_2}{\log w_1}\right)^{\kappa}\left(1+\frac{C}{\log w_1}\right) \quad \text{for} \quad 2\leq w_1\leq w_2.$$

\vspace{1mm}

\noindent
There exist two arithmetic functions $\lambda^{-}$ and $\lambda^{+}$ such that

\begin{itemize}
	\item $\lambda^{\pm}(1)=1, |\lambda^{\pm}|\leq 1,$
	\item $\supp(\lambda^{\pm})\subseteq \{d\mid \prod_{p\in \CP}p:d\leq z^u\},$
	\item $(1\ast\lambda^{-})(n)\leq \mathds{1}_{(\cdot,P(z))=1}(n)\leq (1\ast\lambda^{+})(n)$ for all $n\in\N$ and
\end{itemize}

$$\sum_{d\mid \prod_{p\in\CP}p}\frac{\lambda^{\pm}(d)\nu(d)(\log d)^r}{d}=\!\!\sum_{d\mid \prod_{p\in\CP}p}\frac{\mu(d)\nu(d)(\log d)^r}{d}+O_{r,k}\bigg(\!(\log z)^ru^{-u/2}\prod_{p\in\CP}\left(1-\frac{\nu(p)}{p}\right)\!\!\bigg)\!,$$

\vspace{2mm}

\noindent
for every $r\in\N$.
\end{lemma}

\vspace{1mm}

\begin{proof}
First we show the following claim and then we use it to prove Lemma \ref{sieex}. 

\begin{claim}\label{1}
Let $\nu^*$ be a multiplicative function such that $\nu^*(p)=1$ for $p\!\mid\!P(2ke)$ and $\nu^*(p)=\nu(p)p^{1/\log z}$ when $p\in\CP\cap(2ke,+\infty)$. Then

\begin{equation}\label{last}
\bigg|\sum_{d\mid P(y)}\frac{\mu(d)\nu(md)(\log (md))^r}{d}\bigg|\ll_{r,k}(\log z)^r\nu^*(m)\!\!\prod_{p\mid P(y)}\!\!\left(1-\frac{\nu^*(p)}{p}\right)\!,
\end{equation}

\vspace{2mm}

\noindent
for all $y\in[2,z]$ and for all $m\mid \prod_{p\in\CP,p\geq y}p$.
\end{claim}

\begin{proof}[Proof of Claim \ref{1}]
Let $y\in[2,z]$. For any $s\in\C$ with $|s|\leq 1/\log z$, we may easily show that 

\begin{equation}\label{fors}
\bigg|\prod_{p\in\CP\cap(2ke,y)}\!\!\left(1-\frac{\nu(p)}{p^{1-s}}\right)\!\bigg|\asymp_k\!\!\prod_{p\mid P(y)}\!\!\left(1-\frac{\nu(p)}{p}\right)\!,
\end{equation}

\vspace{2mm}

\noindent
upon noticing that $p^s=1+O(\log p/\log z)$ and that $|\nu(p)/p^{1-s}|\leq 1/2$ for $p\in[2ke,z]\cap\CP.$

We now consider the function $g:\C\rightarrow \C$ with

$$g(s)=m^s\!\!\prod_{p\mid P(y)}\!\!\left(1-\frac{\nu(p)}{p^{1-s}}\right)=\!\sum_{d\mid P(y)}\frac{\mu(d)\nu(d)(md)^s}{d}, \quad s\in\C.$$

\vspace{2mm}

\noindent
For $r\in\N$ and $m\mid \prod_{p\in\CP,p\geq y}p$, Cauchy's residue theorem implies that

\begin{eqnarray}\label{fors2}
\bigg|\sum_{d\mid P(y)}\frac{\mu(d)\nu(md)(\log (md))^r}{d}\bigg|\!\!\!\!\!&=&\!\!\!\!\!\nu(m)|g^{(r)}(0)|=\frac{r!\!\cdot\!\nu(m)}{2\pi}\!\cdot\!\bigg|\int_{|s|=\frac{1}{\log z}}\frac{g(s)}{s^{r+1}}\dee s\bigg|\nonumber\\
&\ll_r&\!\!\!\!(\log z)^r\nu(m)m^{\frac{1}{\log z}}\max\limits_{|s|=1/\log z}\!\bigg|\prod_{p\mid P(y)}\!\!\left(1-\frac{\nu(p)}{p^{1-s}}\right)\!\bigg|\nonumber\\
&\ll_{r,k}&\!\!\!\!(\log z)^r\nu(m)m^{\frac{1}{\log z}}\!\!\prod_{p\mid P(y)}\!\!\left(1-\frac{\nu(p)}{p}\right)\!.
\end{eqnarray}

\vspace{2mm}

\noindent
For the last estimate we first noted that $|\prod_{p\mid P(2ke)}(1-\nu(p)p^{s-1})|\leq \prod_{p\mid P(2ke)}(1+kp^{\Re(s)-1})\leq \prod_{p\mid P(2ke)}(1+kp^{1/\log z-1})\leq \prod_{p\leq 2ke}(1+kp^{1/\log 2-1})=:c_3(k)$ and then we used (\ref{fors}).

Applying (\ref{fors}) with $s=1/\log z$, it follows that

\begin{equation}\label{fors3}
\prod_{p\mid P(y)}\!\!\left(1-\frac{\nu^*(p)}{p}\right)=\!\!\prod_{p\mid P(2ke)}\!\!\bigg(1-\frac{1}{p}\bigg)\prod_{p\in\CP\cap(2ke,y)}\!\!\left(1-\frac{\nu^*(p)}{p}\right)\asymp_k\!\!\prod_{p\mid P(y)}\!\!\left(1-\frac{\nu(p)}{p}\right)\!.
\end{equation}

\vspace{3mm}

\noindent
The term $\prod_{p\mid P(2ke)}(1-1/p)$ was absorbed in the implied constants, because $c_4(k):=\prod_{p\leq 2ke}(1-1/p)\leq \prod_{p\mid P(2ke)}(1-1/p)<1$. We make use of (\ref{fors3}), and so (\ref{fors2}) becomes

\begin{equation}\label{aux}
\bigg|\sum_{d\mid P(y)}\frac{\mu(d)\nu(md)(\log (md))^r}{d}\bigg|\ll_{r,k}(\log z)^r\nu(m)m^{\frac{1}{\log z}}\!\!\prod_{p\mid P(y)}\!\!\left(1-\frac{\nu^*(p)}{p}\right)\!.
\end{equation}

\vspace{2mm}

\noindent
At this point, we write $m=m_1m_2$, where $m_1=\prod_{p\mid m,p\leq 2ke}p$ and $m_2=\prod_{p\mid m,p>2ke}p$. For the positive integer $m_1$ we have 

$$0\leq\nu(m_1)\leq\!\!\prod_{p\mid m,p\leq 2ke}\!\!\nu(p)\leq k^{2ke}\quad \text{and}\quad m_1^{\frac{1}{\log z}}\leq m_1^{\frac{1}{\log 2}}\leq (2ke)^{\frac{2ke}{\log 2}}.$$

\vspace{1mm} 

\noindent
Thus, $\nu(m)m^{1/\log z}\ll_k\nu(m_2)m_2^{1/\log z}=\nu^*(m)$, and combining this inequality with (\ref{aux}), we arrive at the estimate (\ref{last}). So, the claim has now been proven.
\end{proof}

Since $p\leq z$ for any prime $p\in\CP$, if $p\in\CP\cap(2ke,+\infty)$, then $\nu^*(p)\leq ke<p/2$. On the other hand, if $p\in\CP$ and $p\leq 2ke$, then $\nu^*(p)=1<p/2$ as well. So, $0\leq\nu^*(p)<p/2$ and $\nu^*$ is bounded on the primes $p\in\CP$. Now we use the inequalities $0\leq\nu^*(p)<p/2$ for all $p\in\CP$ and  Claim \ref{1} to prove that

$$\sum_{d\mid \prod_{p\in\CP}p}\frac{\lambda^{\pm}(d)\nu(d)(\log d)^r}{d}=\!\!\sum_{d\mid \prod_{p\in\CP}p}\frac{\mu(d)\nu(d)(\log d)^r}{d}+O_{r,k}\bigg(\!(\log z)^ru^{-u/2}\prod_{p\in\CP}\left(1-\frac{\nu^*(p)}{p}\right)\!\!\bigg),$$

\vspace{2mm}

\noindent
where $\lambda^{\pm}$ are arithmetic functions satisfying the conditions mentioned in the statement of the lemma. The proof of the lemma will then be complete, since one may replace $\nu^*$ by $\nu$ in the big-Oh term by using (\ref{fors3}) with $y=z$.

First, let us put $f=\nu\cdot\log^r$ for simplicity. The construction of the functions $\lambda^{\pm}$ is described in \cite[Chapter 19]{dimb}. The proof of the last aforementioned equality is along the same lines as the proof of Theorem 19.1 in \cite{dimb}. One only has to make minor notational changes and rename the $y_j$'s to $z_j$'s and redefine $V(y)$ as 

$$V_f(y)=\sum_{d\mid P(y)}\frac{\mu(d)f(d)}{d}.$$

\noindent
They also have to redefine $V_n$ as

$$V_{n,f}=\sum_{\substack{z_n<p_n<\ldots<p_1\leq z\\p_1,\ldots,p_n\in\CP\\p_i\leq z_i\,(i<n, i\equiv n(\text{mod}\,2))}}\!\!\frac{1}{p_1\cdots p_n}\,\sum_{d\mid P(p_n)}\frac{\mu(d)f(p_1\cdots p_md)}{d}.$$

\vspace{1mm}

\noindent
Then, in an analogous way as in \cite[Theorem 19.1]{dimb}, we have that

\begin{equation}\label{cop1}
V_f(z)-\sum_{d\mid \prod_{p\in\CP}p}\frac{\lambda^+(d)f(d)}{d}=-\sum_{j>J}V_{2j-1,f},
\end{equation}

\vspace{1mm}

\noindent
for the integer $J$ which is defined in \cite[Chapter 19, p.~195, 197]{dimb}. Furthermore, since $z_n<p_n$, using (\ref{last}), we deduce that

\begin{eqnarray}\label{cop2}
|V_n|\!\!\!&\ll_{r,k}&\!\!\!(\log z)^r\!\!\!\sum_{\substack{z_n<p_n<\ldots<p_1\leq z\\p_1,\ldots,p_n\in\CP}}\!\!\frac{\nu^*(p_1)\cdots \nu^*(p_n)}{p_1\cdots p_n}\prod_{p\mid P(p_n)}\!\!\left(1-\frac{\nu^*(p)}{p}\right)\nonumber \\
&\ll_{r,k}&\!\!\!(\log z)^r\!\prod_{p\mid P(z_n)}\!\!\left(1-\frac{\nu^*(p)}{p}\right)\!\sum_{\substack{z_n<p_n<\ldots<p_1\leq z\\p_1,\ldots,p_n\in\CP}}\!\!\frac{\nu^*(p_1)\cdots \nu^*(p_n)}{p_1\cdots p_n}\nonumber \\
&\ll_{r,k}&\!\!\!\frac{(\log z)^r}{n!}\!\cdot\!\prod_{p\mid P(z_n)}\!\!\left(1-\frac{\nu^*(p)}{p}\right)\!\cdot\!\bigg(\sum_{p\in\CP\cap(z_n,z]}\frac{\nu^*(p)}{p}\bigg)^n.
\end{eqnarray}

\vspace{2mm}

\noindent
In the last step we applied the Erd\H{o}s trick to drop the ordering condition from the sum of the second line. Now we can simply use (\ref{cop1}) and (\ref{cop2}) and follow the proof of Theorem 19.1 of \cite{dimb} to show that

$$\sum_{d\mid \prod_{p\in\CP}p}\frac{\lambda^{\pm}(d)\nu(d)(\log d)^r}{d}=\!\!\sum_{d\mid \prod_{p\in\CP}p}\frac{\mu(d)\nu(d)(\log d)^r}{d}+O_{r,k}\bigg(\!(\log z)^ru^{-u/2}\prod_{p\in\CP}\left(1-\frac{\nu^*(p)}{p}\right)\!\!\bigg).$$

\vspace{2mm}

\noindent
The proof for the function $\lambda^-$ is similar. Its only difference is the use of the identity

$$V_f(z)-\sum_{d\mid \prod_{p\in\CP}p}\frac{\lambda^+(d)f(d)}{d}=\sum_{j>J}V_{2j,f}$$

\noindent
in place of (\ref{cop1}).
\end{proof}

\vspace{1mm}

\section{Estimates for Sifted Sums of Multiplicative Functions} 

We begin by proving two preliminary estimates, Lemma \ref{1stl} and Lemma \ref{basle}. These lemmas are essential for the proof of Proposition \ref{bas} which is one of the central results of this section.

\vspace{1mm}

\begin{lemma}\label{1stl}
Let $D\in\N$ and set $k_D=\prod_{p\leq D^3}p$. There exists a real number $c=c(D)\in(0,1)$ such that

\begin{equation}\label{1use}
\sum_{\substack{n\leq x\\(n,k_D)=1}}\!\!\!D^{\Omega(n)}=x\sum_{i=0}^{D-1}a_{i,D}(\log x)^i+O(x^{1-c})
\end{equation}

\vspace{2mm}

\noindent
for some constants $a_{i,D}$ that depend at most on $i$ and $D$.
\end{lemma}

\begin{proof}
This lemma follows from a routine application of the convolution method. We write $\mathds{1}_{(\cdot,k_D)=1}\cdot D^{\Omega}=\tau_D\ast h$ and use the fact that there exists a real number $\eta\in(0,1)$ such that

$$\sum_{n\leq y}\tau_D(n)=y\sum_{i=0}^{D-1}\alpha_{i,D}(\log y)^{i}+O(y^{1-\eta}),$$

\vspace{2mm}

\noindent
for all $y\geq 1$. We leave the details to the reader.
\end{proof}

\vspace{1mm}

\begin{lemma}\label{basle}
	Fix a natural number $D$ and a real number $A>0$. If $x\geq 3$ and $f\in \mathcal{F}(D,A)$, then
	
	\begin{equation}\label{des}
	\sum_{\substack{n\leq x\\(n,d)=1}}\!\!f(n)\ll \frac{x}{(\log x)^A}\left(\frac{d}{\phi(d)}\right)^D\!\! \quad \text{for} \quad d\leq x.
	\end{equation}
	
\end{lemma}

\begin{proof}
	For any $n$, there is a unique way to write it as $n=am$, where all prime divisors of $a$ divide $d$ and $(m,d)=1$. Consequently, we find that 
	
	$$S(x):=\sum_{n\leq x}f(n)=\sum_{n\leq x}h(n)\!\!\sum_{\substack{m\leq x/n\\(m,d)=1}}\!\!f(m),$$
	
	\noindent
	where $h$ is the multiplicative function with $h(p^{\nu})=f(p^{\nu})\mathds{1}_{p\mid d}$ for $\nu\in \mathbb{N}$. By applying Lemma \ref{mobin}, we get that
	
	\begin{equation}\label{foru}
	\sum_{\substack{n\leq x\\(n,d)=1}}\!\!\!f(n)=\sum_{n\leq x}h^{-1}(n)S(x/n),\,\,\,
	\end{equation}
	
	\noindent
	where $h^{-1}$ denotes the Dirichlet inverse of $h$. Note that $h\in \mathcal{F}(D)$. Therefore, $h^{-1}\in \mathcal{F}(D)$ too. Furthermore, one observes that $h^{-1}(p^{\nu})=f^{-1}(p^{\nu})\mathds{1}_{p\mid d}$ for any $\nu\in\mathbb{N}$. 
	
	We now split the sum of the right-hand side of (\ref{foru}) into the two parts
	
	$$T_1=\sum_{n\leq \sqrt{x}}\!h^{-1}(n)S(x/n)\! \quad \text{and} \quad T_2=\!\!\!\sum_{\sqrt{x}<n\leq x}\!\!\!h^{-1}(n)S(x/n).$$
	
	We begin with the estimation of $T_1$. Since $h^{-1}\in \mathcal{F}(D)$, it is true that $|h^{-1}|\leq {\tau}_D$, and so
	
	\begin{eqnarray}\label{T1}
	\,|T_1|\leq\sum_{\substack{n\leq \sqrt{x}\\p\mid n \Rightarrow p\mid d}}\!\!\!\tau_D(n)|S(x/n)|\!\!\!&\ll&\!\!\!\frac{x}{(\log x)^A}\sum_{\substack{n\leq \sqrt{x}\\p\mid n \Rightarrow p\mid d}}\!\!\!\frac{\tau_D(n)}{n} \nonumber \\
	&\leq& \!\!\!\frac{x}{(\log x)^A}\prod_{p\mid d}\left(1+\frac{\tau_D(p)}{p}+\frac{\tau_D(p^2)}{p^2}+\ldots\right) \nonumber \\
	&=&\!\!\!\frac{x}{(\log x)^A}\prod_{p\mid d}\left(\sum_{j\geq 0}\binom{D+j-1}{j}\frac{1}{p^{\,j}}\right)\!=\frac{x}{(\log x)^A}\left(\frac{d}{\phi(d)}\right)^D\!\!.
	\end{eqnarray}
	
	\vspace{3mm}
	
	We continue by bounding the sum $T_2$. Since $\tau_D(n)\ll n^{1/4}\leq x^{1/4}$ for $n\leq x$, we have
	
	\begin{eqnarray*}
		|T_2|\!\!&\leq& \!\!\!\!\sum_{\substack{\sqrt{x}<n\leq x\\p\mid n \Rightarrow p\mid d}}\!\!\!\tau_D(n)|S(x/n)|\ll x\!\!\sum_{\substack{\sqrt{x}<n\leq x\\p\mid n \Rightarrow p\mid d}}\!\!\!\frac{\tau_D(n)}{n}<\sqrt{x}\!\!\sum_{\substack{\sqrt{x}<n\leq x\\p\mid n \Rightarrow p\mid d}}\!\!\!\tau_D(n)\ll x^{3/4}\!\!\sum_{\substack{n\leq x\\p\mid n \Rightarrow p\mid d}}\!\!1.
	\end{eqnarray*}
	
	\noindent
	Proposition \ref{geo} implies that
	
	$$\sum_{\substack{n\leq x\\p\mid n\Rightarrow p\mid d}}\!\!1\leq \omega(d)\cdot\frac{(\log x+\log d)^{\omega(d)}}{\omega(d)!\left(\prod_{p\mid d}\log p\right)}<4^{\omega(d)}\frac{(\log x)^{\omega(d)}}{\log 2\cdot\omega(d)!}\ll 32^{\omega(d)}x^{1/8},$$
	
	\vspace{1mm}
	
	\noindent
	and so $T_2\ll 32^{\omega(d)}x^{7/8}$. Using the maximal order of $\omega$, we have that $\omega(d)\ll \log x/\log \log x$ for $d\leq x$. Therefore, there exists a constant $C>0$ such that
	
	\begin{equation}\label{T2}
	T_2\ll x^{\frac{7}{8}+\frac{C}{\log\log x}}\ll \frac{x}{(\log x)^A}\left(\frac{d}{\phi(d)}\right)^D\!\!.
	\end{equation}
	
	\vspace{4mm}
	
	Combination of (\ref{foru}), (\ref{T1}) and (\ref{T2}) completes the proof of the lemma.
\end{proof}

\vspace{1mm}

Now that Lemmas \ref{1stl} and \ref{basle} are proven, we combine them with Lemma \ref{sieex} to establish an upper bound for the sifted partial sums $\sum_{n\leq x,P^-(n)>z}f(n)$ of a function $f\in\CF(D,A)$, where $D\in\N$ and $A>0$.

\vspace{1mm}

\begin{proposition}\label{bas}
	Fix a natural number $D$ and a positive real number $A$. If $x\geq 3$ and $f\in \mathcal{F}(D,A)$, then there exists a constant $\alpha=\alpha(D)\in(0,1)$ such that
	
	$$\sum_{\substack{n\leq x\\P^-(n)>z}}\!\!\!f(n)\ll \frac{x(\log z)^D}{(\log x)^A}+\frac{x^{1-\alpha/\log z}}{\log z}\,\,\,$$
	
	\noindent
	for all $z\in[2,x]$.
\end{proposition}

\begin{proof}
	Let $C=\min\{1/16,c/2\}$, where $c$ is the constant appearing in Lemma \ref{1stl}. First we show that the estimate holds trivially when $z>x^{C/(4D+1)}$. Indeed, in this case $\log x/\log z\asymp 1$, and so using Theorem \ref{fromb} with the divisor-bounded, multiplicative function $\tau_D\cdot \mathds{1}_{P^-(\cdot)>z}$, we conclude that
	
	$$\Bigg|\sum_{\substack{n\leq x\\P^-(n)>z}}\!\!\!f(n)\Bigg|\leq \!\sum_{\substack{n\leq x\\P^-(n)>z}}\!\!\!\tau_D(n)\ll \frac{x}{\log x}\left(\frac{\log x}{\log z}\right)^D\!\!\ll\frac{x^{1-C/(2\log z)}}{\log z}.$$ 
	
	\vspace{2mm}
	
	Now it remains to prove the proposition when $z\leq x^{C/(4D+1)}$. Assuming that $x$ is large enough in terms of $D$, when $z\leq D^3$, we can use Lemma \ref{basle}. For $z>D^3$, the condition $P^-(n)>z$ implies that $(n,k_D)=1$, where $k_D=\prod_{p\leq D^3}p$. So, using the arithmetic functions $\lambda^{-}$ and $\lambda^{+}$ of Lemma \ref{sieex} with $u=C\log x/\log z$, we write
	
	\begin{eqnarray}\label{sie}
	\,\,\,\,\,\,\,\,\,\,\,\,\,\,\,\,\,\,\,\,\,\,\,\,\sum_{\substack{n\leq x\\P^-(n)>z}}\!\!\!f(n)=\!\!\!\sum_{\substack{n\leq x\\(n,k_D)=1}}\!\!(\lambda^+\ast 1)(n)f(n)+O\Big(\!\sum_{\substack{n\leq x\\(n,k_D)=1}}\!\!(\lambda^+\ast 1-\lambda^-\ast 1)(n)|f(n)|\,\Big).
	\end{eqnarray}
	
	\vspace{2mm}
	
	Since $|f(n)|\leq \tau_D(n)\leq D^{\Omega(n)}$ for all $n\in\mathbb{N}$, it follows that
	
	\begin{equation}\label{opcon}
	\sum_{\substack{n\leq x\\(n,k_D)=1}}\!\!(\lambda^+\ast 1-\lambda^-\ast 1)(n)|f(n)|\,\leq\!\sum_{(d,k_D)=1}\!\!(\lambda^+(d)-\lambda^-(d))D^{\Omega(d)}\!\!\!\sum_{\substack{m\leq x/d\\(m,k_D)=1}}\!\!\!D^{\Omega(m)}.
	\end{equation}
	
	\noindent
    We insert the formula of Lemma \ref{1stl} in the right-hand side of (\ref{opcon}) and use the binomial theorem to expand the resulting powers $(\log (x/d))^{\nu}=(\log x-\log d)^{\nu}$ with $\nu\leq D-1$. Since $|\lambda^{\pm}|\leq 1$, the contribution coming from the error term of (\ref{1use}) is
	
	$$\ll x^{1-c}\sum_{d\leq x^C}D^{\Omega(d)}\ll x^{1-c+C}(\log x)^{D-1}\leq x^{1-c/2}(\log x)^{D-1}\ll \frac{x(\log z)^D}{(\log x)^A}.$$
	
	\vspace{2mm}
	
	\noindent
	The summands coming from the main term of (\ref{1use}) contain expressions of the form
	
	$$\sum_{d}\frac{(\lambda^+(d)-\lambda^-(d))D^{\Omega(d)}\mathds{1}_{(d,k_D)=1}(\log d)^j}{d}$$
	
	\vspace{2mm}
	
	\noindent
	for $j\in \{0,\ldots,D-1\}$. Each one of these expressions is multiplied by a logarithmic factor $(\log x)^{\ell}$ with $\ell+j\leq D-1$. Since Lemma \ref{sieex} implies that
	
	$$\sum_{(d,k_D)=1}\!\!\frac{(\lambda^+(d)-\lambda^-(d))D^{\Omega(d)}(\log d)^j}{d}\ll \frac{x^{-\frac{C}{\log z}}}{{(\log z)^{D-j}}} \quad \text{for} \quad \!j\in\{0,\ldots,D-1\},$$
	
	\vspace{2mm}
	
	\noindent
	we finally infer that
	
	$$\sum_{n\leq x}\,(\lambda^+\ast 1-\lambda^-\ast 1)(n)|f(n)|\ll\frac{x^{1-\frac{C}{\log z}}(\log x)^{D-1}}{(\log z)^{D}}\ll \frac{x^{1-\frac{C}{2\log z}}}{\log z}.$$
	
	\vspace{1mm}
	
	\noindent
	So, if we define $\alpha:=C/2$,  relation (\ref{sie}) becomes
	
	\begin{equation}\label{n4.6}
	\sum_{\substack{n\leq x\\P^-(n)>z}}\!\!\!f(n)=\!\!\sum_{\substack{n\leq x\\(n,k_D)=1}}\!\!(\lambda^+\ast 1)(n)f(n)+O\left(\frac{x(\log z)^D}{(\log x)^A}+\frac{x^{1-\alpha/\log z}}{\log z}\right)\!.\,\,\,
	\end{equation}
	
	\vspace{2mm}
	
	We now turn to the estimation of the ``main term'' of (\ref{n4.6}). For every two natural numbers $m,d$, there is a unique way to write $md=m'd'$ where $(m',d')=1, d\mid d'$ and all the prime factors of $d'$ divide $d$. Then
	
	$$\sum_{\substack{n\leq x\\(n,k_D)=1}}\!\!(\lambda^+\ast 1)(n)f(n)=\!\!\!\sum_{\substack{md\leq x\\(md,k_D)=1}}\!\!\!\!\lambda^+(d)f(md)=\!\!\sum_{(d,k_D)=1}\!\!\lambda^+(d)\!\sum_{\substack{d'\leq x,d\mid d'\\p\mid d\Leftrightarrow p\mid d'\\(k_D,d')=1}}\!\!\!f(d')\!\!\!\sum_{\substack{m'\leq x/d'\\(dk_D,m')=1}}\!\!\!\!f(m').$$
	
	\noindent
	We divide the inner sum on $d'$ into two sums $S_1$ and $S_2$. In $S_1$ we are summing over the range $d'\leq \sqrt{x}$. In the sum $S_2$ we have $\sqrt{x}<d'\leq x$. We apply Lemma \ref{basle} with $x$ large enough to the sums 
	
	$$\sum_{\substack{m'\leq x/d'\\(m',dk_D)=1}}\!\!\!\!f(m').$$
	
	\vspace{2mm}
	
	\noindent
	Then, the sum on $d$, coming from $S_1$, is 
	
	\begin{equation}\label{es1}
	\ll \frac{x}{(\log x)^A}\sum_d\left(\frac{d}{\phi(d)}\right)^D\!\!\!\sum_{\substack{d'\!,\,d\mid d'\\p\mid d\Leftrightarrow p\mid d'}}\!\!\!\frac{\tau_D(d')}{d'}
	\end{equation}
	
	\noindent
	and the sum on $d$, coming from $S_2$, is
	
	\begin{equation}\label{es2}
	\ll\sqrt{x}\,\sum_d\left(\frac{d}{\phi(d)}\right)^D\!\!\!\sum_{\substack{d'\leq x\\p\mid d\Leftrightarrow p\mid d'}}\!\!\!\tau_D(d').
	\end{equation}
	
	\noindent
	So, in order for the proof to be completed, we need to show that the quantities of (\ref{es1}) and (\ref{es2}) are $\ll x(\log z)^D/(\log x)^A$. 
	
	First, we have that
	
	$$\sum_{\substack{d'\!,\,d\mid d'\\p\mid d\Leftrightarrow p\mid d'}}\!\!\!\frac{\tau_D(d')}{d'}\leq\prod_{p\mid d}\left(\sum_{j\geq 1}\frac{\tau_D(p^j)}{p^j}\right)=\prod_{p\mid d}\left(\left(1-\frac{1}{p}\right)^{-D}\!\!-1\right)
	\leq \frac{D^{\omega(d)}}{d}\left(\frac{d}{\phi(d)}\right)^{D+1}\!\!,$$
	
	\vspace{2mm}
	
	\noindent
	where the last estimate follows from the inequality $(1-1/p)^{-D}-1\leq Dp^{-1}(1-1/p)^{-D-1}$, which may be obtained by applying the Mean Value Theorem to the function $t\mapsto t^{-D}$. Therefore, the expression of (\ref{es1}) is bounded by 
	
	\begin{eqnarray*}
		\frac{x}{(\log x)^A}\sum_{d\mid \prod_{p\leq z}p}\frac{D^{\omega(d)}}{d}\left(\frac{d}{\phi(d)}\right)^{2D+1}\!\!\!\!&\leq& \!\!\!\frac{x}{(\log x)^A}\prod_{p\leq z}\left(1+\frac{D}{p}\left(1+\frac{1}{p-1}\right)^{2D+1}\right)\\
		&=&\!\!\!\frac{x}{(\log x)^A}\prod_{p\leq z}\left(1+\frac{D}{p}+O\left(\frac{1}{p^2}\right)\right)\\
		&\ll& \!\!\!\frac{x}{(\log x)^A}\exp\left(D\sum_{p\leq z}\frac{1}{p}\right)\ll \frac{x(\log z)^D}{(\log x)^A}.
	\end{eqnarray*}

\vspace{3mm}
	
	We now continue with the estimation of the expression of (\ref{es2}). As in the proof of Lemma \ref{basle}, Proposition \ref{geo} implies that
	
	$$\sum_{\substack{d'\\p\mid d\Leftrightarrow p\mid d'}}\!\!1<\,4^{\omega(d)}\frac{(\log x)^{\omega(d)}}{\log 2\cdot\omega(d)!}\ll 16^{\omega(d)}x^{1/4}.$$
	
	\vspace{1mm}
	
	\noindent
	We also have the inequalities $d/\phi(d)\ll \log\log d\leq \log\log x$ and $\tau_D(d)\ll d^{1/8}\leq x^{1/8}$ for $d\leq x$. Consequently, the expression of (\ref{es2}) is
	
	\begin{eqnarray*}
		\ll x^{7/8}(\log \log x)^D\sum_{d\leq x^C}16^{\omega(d)}\ll x^{15/16}(\log x)^{15}(\log\log x)^D\ll\frac{x(\log z)^D}{(\log x)^A}
	\end{eqnarray*}

\vspace{2mm}
	
    \noindent
	and the proposition is finally proved.
\end{proof}

\vspace{1mm}

The next lemma is a rather technical result and is useful for the proof of Proposition \ref{corpr}. It concerns a bound for the tails of a Dirichlet series and its proof is an easy application of partial summation. 

\vspace{1mm}

\begin{lemma}\label{tails}
Fix a natural number $D$ and two real numbers $\varepsilon>0$ and $z\geq 2$. Let $\Gamma$ be a mutliset of $m$ elements, counting the multiplicities. Let also $f$ be an arithmetic function and suppose that there exist some $\delta\in(0,1)$ and some $A\geq m+1+\varepsilon$ such that

\begin{equation}\label{hyp}
\sum_{\substack{n\leq x\\P^-(n)>z}}\!\!\!f_{\Gamma}(n)\ll\frac{x(\log z)^{D-m}}{(\log x)^{A-m}}+\frac{x^{1-\delta/\log z}}{\log z}\quad \text {whenever}\quad x\geq z,
\end{equation}

\noindent
where $f_{\Gamma}=f\ast \tau_{\Gamma}$. For $N\geq \max\{3,z\}$ and $s=\sigma+it$ with $\sigma\in[1,2]$ and $t\in\R$, we have

\begin{equation*}
\sum_{\substack{n>N\\P^-(n)>z}}\!\!\!\frac{f_{\Gamma}(n)}{n^s}\ll_{\varepsilon, \delta} (1+|t|)N^{1-\sigma}\left(\frac{(\log z)^{D-m}}{(\log N)^{A-m-1}}+N^{-\frac{\delta}{\log z}}\right)\!.
\end{equation*}
\end{lemma}

\begin{proof}
Let $M>N$. Then, partial summation implies that

\begin{eqnarray}\label{part}
\,\,\,\,\,\,\,\,\,\,\,\,\,\,\sum_{\substack{N<n\leq M\\P^-(n)>z}}\!\!\frac{f_{\Gamma}(n)}{n^s}=\Big(\!\!\sum_{\substack{n\leq x\\P^-(n)>z}}\!\!\!f_{\Gamma}(n)\Big)x^{-s}\Big|_{x=N}^{M}+s\int_N^M\Big(\!\!\sum_{\substack{n\leq y\\P^-(n)>z}}\!\!\!f_{\Gamma}(n)\Big)\frac{\dee y}{y^{s+1}}.
\end{eqnarray}

\vspace{1mm}

Using the hypothesis (\ref{hyp}) twice, once with $x=N$ and once with $x=M$, we conclude that

\begin{eqnarray}\label{1stp}
\,\,\,\,\,\,\,\,\,\Big(\!\!\sum_{\substack{n\leq x\\P^-(n)>z}}\!\!\!f_{\Gamma}(n)\Big)x^{-s}\Big|_{x=N}^{M}\ll\frac{(\log z)^{D-m}}{(\log N)^{A-m}N^{\sigma-1}}+\frac{N^{-\frac{\delta}{\log z}}}{N^{\sigma-1}(\log z)}.
\end{eqnarray}

We now focus on the integral of the right-hand side of (\ref{part}). The hypothesis (\ref{hyp}) yields

\begin{eqnarray*}
\,\,\,\,\,\,\,\,\,\,\,\int_N^M\Big(\!\!\sum_{\substack{n\leq y\\P^-(n)>z}}\!\!\!f_{\Gamma}(n)\Big)\frac{\dee y}{y^{s+1}}\ll (1+|t|)\Big(\int_N^M\frac{(\log z)^{D-m}}{(\log y)^{A-m}y^{\sigma}}\dee y+\frac{1}{\log z}\int_N^M\frac{\dee y}{y^{\frac{\delta}{\log z}+\sigma}}\Big).
\end{eqnarray*}

\vspace{1mm}

\noindent
We also have that

\begin{equation*}
\int_N^M\frac{(\log z)^{D-m}}{(\log y)^{A-m}y^{\sigma}}\dee y\leq \frac{(\log z)^{D-m}}{N^{\sigma-1}}\int_N^{\infty}\frac{\dee y}{(\log y)^{A-m}y}\ll_{\varepsilon} \frac{(\log z)^{D-m}}{(\log N)^{A-m-1}N^{\sigma-1}}.\,\,\,
\end{equation*}

\vspace{3mm}

\noindent
Furthermore,

\begin{eqnarray*}
\,\,\,\,\int_N^M\frac{\dee y}{y^{\frac{\delta}{\log z}+\sigma}}\leq \int_N^{\infty}\frac{\dee y}{y^{\frac{\delta}{\log z}+\sigma}}=\frac{N^{-\frac{\delta}{\log z}}}{N^{\sigma-1}(\sigma+\delta/\log z-1)}\leq \delta^{-1}(\log z)N^{1-\sigma-\frac{\delta}{\log z}}.
\end{eqnarray*}

\vspace{4mm}

\noindent
Therefore, we infer that

\begin{equation}\label{ins2}
\int_N^M\Big(\!\!\sum_{\substack{n\leq y\\P^-(n)>z}}\!\!\!f_{\Gamma}(n)\Big)\frac{\dee y}{y^{s+1}}\ll_{\varepsilon,\delta}(1+|t|)N^{1-\sigma}\left(\frac{(\log z)^{D-m}}{(\log N)^{A-m-1}}+N^{-\frac{\delta}{\log z}}\right)\!.\,\,\,\,\,\,\,\,\,
\end{equation}

We insert the estimates (\ref{1stp}) and (\ref{ins2}) in (\ref{part}) and obtain the desired inequality. 
\end{proof}

\vspace{1mm}

We close this section with the following stronger version of Proposition \ref{bas} which will serve as the first basic ingredient going into the proof of Theorem \ref{mthm}.

\vspace{1mm} 

\begin{proposition}\label{corpr}
Suppose that $D\in\N$ and that $A>D+1$. Let also $f\in \CF(D,A)$ and $\tilde{\Gamma}$ be a (multi-)subset of the multiset $\Gamma$ of the ordinates of the zeros of $L(\cdot,f)$ on $\Re(s)=1$. Using the definition of $V_t$ for $t\in\R$ before Lemma \ref{dle}, if $x\geq z\geq V_{\tilde{\Gamma}}:=\max_{\gamma\in\tilde{\Gamma}}V_{\gamma}$, then there exists a real number $\kappa=\kappa(D)\in(0,1)$ such that
	
\begin{equation*}
\sum_{\substack{n\leq x\\P^-(n)>z}}\!\!\!f_{\tilde{\Gamma}}(n)\ll_{\tilde{\Gamma}} \frac{x(\log z)^{D-m}}{(\log x)^{A-m}}+\frac{x^{1-\kappa/\log z}}{\log z},
\end{equation*}

\vspace{1mm}
	
\noindent
where $m$ is the number of elements of $\,\tilde{\Gamma}$ with the multiplicities being counted.
\end{proposition}

\begin{proof}
We perform induction on the number of elements $m$ of a multisubset of $\Gamma$. The proposition holds when $m=0$ because of Proposition \ref{bas}. For $m\geq 1$, we assume that the proposition is true for any multisubset of $\Gamma$ with $m-1$ elements. We will show that it remains true for a multisubset $\tilde{\Gamma}$ of $m$ elements.

When $\sqrt{x}<z$, then $\log x\asymp\log z$ and we may argue as in the beginning of the proof of Proposition \ref{bas} to show that

$$\sum_{\substack{n\leq x\\P^-(n)>z}}\!\!\!f_{\tilde{\Gamma}}(n)\ll_{\,\tilde{\Gamma},\varepsilon} \frac{x^{1-\varepsilon/\log z}}{\log z},$$

\noindent
for any $\varepsilon>0$. 

If $\sqrt{x}\geq z$, let $\gamma$ be an element of $\tilde{\Gamma}$ and write $\tilde{\Gamma}=\Gamma'\cup\{\gamma\}$. We have that the cardinality of $\Gamma'$ is $m-1$ and that $f_{\tilde{\Gamma}}(n)=\sum_{ab=n}f_{\Gamma'}(a)b^{i\gamma}$ for all $n\in\N$. So,

\begin{equation}
\sum_{\substack{n\leq x\\P^-(n)>z}}\!\!\!f_{\tilde{\Gamma}}(n)=\!\!\sum_{\substack{a\leq \sqrt{x}\\P^-(a)>z}}\!\!\!f_{\Gamma'}(a)\!\!\sum_{\substack{b\leq x/a\\P^-(b)>z}}\!\!\!b^{i\gamma}\,+\!\sum_{\substack{b\leq \sqrt{x}\\P^-(b)>z}}\!\!\!b^{i\gamma}\!\sum_{\substack{\sqrt{x}<a\leq x/b\\P^-(a)>z}}\!\!\!f_{\Gamma'}(a):=S_1+S_2,
\end{equation}

\vspace{1mm}

\noindent
say. Since $\sqrt{x}\geq z$, Lemma \ref{dle} gives that

\begin{eqnarray}\label{com1}
\,\,\,\,\,\,\,\,\,\,\,\,\,\,\,\,\,\,\,\,\,\,\,\,\,\,\,\,\,\,\,S_1=\frac{x^{1+i\gamma}}{1+i\gamma}\,\prod_{p\leq z}\bigg(1-\frac{1}{p}\bigg)\!\!\!\sum_{\substack{a\leq \sqrt{x}\\P^-(a)>z}}\!\!\!\frac{f_{\Gamma'}(a)}{a^{1+i\gamma}}+O\bigg(\frac{x^{1-\frac{1}{30\log z}}}{\log z}\!\!\sum_{\substack{a\leq \sqrt{x}\\P^-(a)>z}}\!\!\frac{\tau_{D+m-1}(a)}{a^{1-\frac{1}{30\log z}}}\bigg).
\end{eqnarray}

\vspace{1mm}

\noindent
First we bound the sum in the big-Oh term. From Mertens's theorem it follows that

\begin{eqnarray}\label{com2}
\sum_{\substack{a\leq \sqrt{x}\\P^-(a)>z}}\!\!\frac{\tau_{D+m-1}(a)}{a^{1-\frac{1}{30\log z}}}\!\!\!&\leq&\!\!\!x^{\frac{1}{60\log z}}\!\!\sum_{\substack{a\leq \sqrt{x}\\P^-(a)>z}}\!\!\frac{\tau_{D+m-1}(a)}{a}\leq x^{\frac{1}{60\log z}}\!\!\prod_{z<p\leq \sqrt{x}}\left(\sum _{j\geq 0}\frac{\tau_{D+m-1}(p^{\,j})}{p^{\,j}}\right)\nonumber\\&=&\!\!\!x^{\frac{1}{60\log z}}\!\!\prod_{z<p\leq \sqrt{x}}\left(1-\frac{1}{p}\right)^{-(D+m-1)}\!\!\ll x^{\frac{1}{60\log z}}\left(\frac{\log x}{\log z}\right)^{D+m-1}\ll_{\tilde{\Gamma}}x^{\frac{1}{40\log z}}.
\end{eqnarray}

\vspace{1mm}

\noindent
We continue by bounding the sum outside the big-Oh term on the right-hand side of (\ref{com1}). Since we have that $V_{\Gamma'}\leq V_{\tilde{\Gamma}}\leq z$ and that $\Gamma'$ contains $m-1$ elements, the inductive hypothesis implies that for all $w\geq z$ it is true that

\begin{equation}\label{com4}
\sum_{\substack{a\leq w\\P^-(a)>z}}\!\!\!f_{\Gamma'}(a)\ll_{\Gamma'}\frac{w\,(\log z)^{D-m+1}}{(\log w)^{A-m+1}}+\frac{w^{1\,-\,\kappa_0/\log z}}{\log z},
\end{equation}

\noindent
where $\kappa_0=\kappa_0(D)$ is some real number of $(0,1)$. We have that $m\leq D$ because of part (b) of Proposition \ref{dsp}. So, $A>(m-1)+2$. Therefore, we can use (\ref{com4}) to apply Lemma \ref{tails} with $\varepsilon=1, \delta=\kappa_0, N=\sqrt{x}$ and $\sigma=1$ and deduce that

$$\,\,\sum_{\substack{a>\sqrt{x}\\P^-(a)>z}}\!\!\!\frac{f_{\Gamma'}(a)}{a^{1+i\gamma}}\ll_{\tilde{\Gamma}}\frac{(\log z)^{D-m+1}}{(\log x)^{A-m}}+x^{-\kappa_0/(2\log z)}.$$

\noindent
Since $\gamma$ is an element of $\tilde{\Gamma}$, the complex number $1+i\gamma$ is a zero of $L(\cdot,f)$. In addition, the multiplicity of $\gamma$ in $\Gamma'$ is smaller than its multiplicity in $\Gamma$, because $\tilde{\Gamma}=\Gamma'\cup\{\gamma\}$ and $\tilde{\Gamma}$ is a multisubset of $\Gamma$. Therefore, $L(1+i\gamma,f_{\Gamma'})=0$, and so

\begin{eqnarray}\label{com5}
\sum_{\substack{a\leq\sqrt{x}\\P^-(a)>z}}\!\!\!\frac{f_{\Gamma'}(a)}{a^{1+i\gamma}}\!\!\!&=&\!\!\!L(1+i\gamma,f_{\Gamma'})\cdot\prod_{p\leq z}\bigg(\sum_{j\geq 0}\frac{f_{\Gamma'}(p^{\,j})}{p^{\,j(1+i\gamma)}}\bigg)^{\!-1}\!-\!\!\sum_{\substack{a>\sqrt{x}\\P^-(a)>z}}\!\!\!\frac{f_{\Gamma'}(a)}{a^{1+i\gamma}}\nonumber\\
&=&\!\!\!-\!\!\sum_{\substack{a>\sqrt{x}\\P^-(a)>z}}\!\!\!\frac{f_{\Gamma'}(a)}{a^{1+i\gamma}}\ll_{\tilde{\Gamma}}\frac{(\log z)^{D-m+1}}{(\log x)^{A-m}}+x^{-\kappa_0/(2\log z)}.
\end{eqnarray}

\noindent
Now, combination of (\ref{com1}) with the estimates (\ref{com2}) and (\ref{com5}) gives

$$S_1\ll_{\tilde{\Gamma}}\frac{x(\log z)^{D-m}}{(\log x)^{A-m}}+\frac{x^{1-\kappa_1/\log z}}{\log z},$$

\noindent
for some $\kappa_1=\kappa_1(D)\in(0,1)$.

It only remains to estimate

$$S_2=\sum_{\substack{b\leq \sqrt{x}\\P^-(b)>z}}\!\!\!b^{i\gamma}\!\sum_{\substack{\sqrt{x}<a\leq x/b\\P^-(a)>z}}\!\!\!f_{\Gamma'}(a).$$

\vspace{1mm}

\noindent
In fact, we are going to show that $S_2$ satisfies the same bound as $S_1$ and then the proof of the proposition will be complete. In the innermost sum of $S_2$, we have that $x/b\geq \sqrt{x}\geq z$. So, using the inductive hypothesis once for $\sqrt{x}$ and once for $x/b$, we get that

\begin{equation}\label{lafin}
S_2\ll_{\Gamma'}\frac{x(\log z)^{D-m+1}}{(\log x)^{A-m+1}}\!\!\sum_{\substack{b\leq \sqrt{x}\\P^-(b)>z}}\frac{1}{b}+\frac{x^{1-\kappa_0/\log z}}{\log z}\!\!\sum_{\substack{b\leq \sqrt{x}\\P^-(b)>z}}\frac{1}{b^{1-\frac{\kappa_0}{\log z}}}.
\end{equation}

\noindent
For $u\geq z$, it is known that $\#\{n\leq u:P^-(n)>z\}\!\ll\!u/\log z$. This estimate and partial summation imply that $\sum_{b\leq \sqrt{x}, P^-(b)>z}1/b\ll \log x/\log z$. Then,

$$\sum_{\substack{b\leq \sqrt{x}\\P^-(b)>z}}\frac{1}{b^{1-\frac{\kappa_0}{\log z}}}\leq x^{\,\kappa_0/(2\log z)}\!\!\sum_{\substack{b\leq \sqrt{x}\\P^-(b)>z}}\frac{1}{b}\ll x^{\,\kappa_0/(2\log z)}\frac{\log x}{\log z}\ll x^{\,2\kappa_0/(3\log z)}.$$

\vspace{1mm}

\noindent
So, finally, the estimate (\ref{lafin}) becomes

$$S_2\ll_{\tilde{\Gamma}}\frac{x(\log z)^{D-m}}{(\log x)^{A-m}}+\frac{x^{1-\kappa_0/(3\log z)}}{\log z}$$

\vspace{2mm}

\noindent
and the proof is finished with $\kappa:=\min\{\kappa_0/3,\kappa_1\}$.
\end{proof}

\vspace{1mm}

\section{Proof of Theorem \ref{mthm}}

We now move on to the proof of Theorem \ref{mthm}. First, we introduce some auxiliary notation. If $g$ is an arithmetic function and $x\geq 2$ is a real number, we define

$$S(x,g)=\sum_{n\leq x}g(n).$$

\vspace{1mm}

\noindent
Moreover, for $z\geq 1$, we set $g_z(n)=g(n)$ when $P^-(n)>z$ and $g_z(n)=0$ otherwise.\\

	The multiset $\Gamma$ in the statement of the theorem consists of the ordinates of the zeros of $L(\cdot,f)$ on the vertical line $\sigma=1$. If $L(\cdot,f)$ has a single root $1+i\gamma$ of multiplicity $D$, then Theorem \ref{mthm} follows directly from Theorem \ref{1r}. So, for the rest of this proof we assume that the largest multiplicity of the roots of $L(\cdot,f)$ on $\Re(s)=1$ is at most $D-1$. Let $m$ be the number of elements of $\Gamma$, with the multiplicities being counted. In Proposition \ref{corpr} we take $\tilde{\Gamma}=\Gamma$ and
	
	\begin{eqnarray*}
	z=x^{1/(L\log u)}\!\!\!&\text{with}&\!\! L=8(\min\{\kappa,\log \log 3\cdot(A-D-1)\})^{-1}\,\, \text{and}\\
	u\!\!\!\!\!\!\!&=&\!\!\!\!\!\!\!\min\{(\log x)^{A-D-1},\sqrt{T}\}.
	\end{eqnarray*}

    \vspace{2mm}
	
	\noindent
	We further assume that $x$ is large enough in terms of $\Gamma$. Indeed, for bounded $x$, Theorem \ref{mthm} holds trivially by adjusting the implied constant in its statement. Then, for $w\in[x^{1/4},x]\subseteq [z^{1/4},x]$, we have
	
	\begin{equation}\label{app}
	\sum_{\substack{n\leq w\\P^-(n)>z}}\!\!\!f_{\Gamma}(n)\ll_{\Gamma} \frac{w}{u\log x}.
	\end{equation}
	
	\vspace{2mm}
	
	\noindent
	Consequently,
	
	$$\sum_{\substack{n\leq x\\P^-(n)>z}}\!\!\!f_{\Gamma}(n)\log n=O(x^{1/3})+\!\!\sum_{\substack{x^{1/4}<n\leq x\\P^-(n)>z}}\!\!\!f_{\Gamma}(n)\log n\ll_{\Gamma} \frac{x}{u},$$
	
	\noindent
	where we used the inequality $|f_{\Gamma}|\leq \tau_{D+m}$ to control the summands with $n\leq x^{1/4}$, whereas we used partial summation and (\ref{app}) for the summands with $n>x^{1/4}$. Now, combining this bound with Dirichlet's hyperbola method applied to the convolution $(f_{\Gamma})_z\ast \Lambda_{(f_{\Gamma})_z}$, we get
	
	\begin{eqnarray}\label{app2}
	\,\,\,\,\sum_{\substack{n\leq x^{1/4}\\P^-(n)>z}}\!\!\!f_{\Gamma}(n)S(x/n,\Lambda_{(f_{\Gamma})_z})\!\!\!&=&\!\!\!S(x^{1/4},(f_{\Gamma})_z)S(x^{3/4},\Lambda_{(f_{\Gamma})_z})\\&-&\!\!\!\sum_{n\leq x^{3/4}}\!\!\Lambda_{(f_{\Gamma})_z}(n)S(x/n,(f_{\Gamma})_z)\nonumber+O_{\Gamma}\left(\frac{x}{u}\right).
	\end{eqnarray}
	
	\noindent
	Note that $\Lambda_{(f_{\Gamma})_z}=(\Lambda_{f_{\Gamma}})_z$. So, since $|\Lambda_{f_{\Gamma}}|\leq (D+m)\cdot\Lambda$, we deduce that $|S(x^{3/4},\Lambda_{(f_{\Gamma})_z})|\leq \sum_{n\leq x^{3/4}}|\Lambda_{f_{\Gamma}}(n)|\mathds{1}_{P^-(n)>z}\leq \sum_{n\leq x^{3/4}}|\Lambda_{f_{\Gamma}}(n)|\leq (D+m)\sum_{n\leq x^{3/4}}\Lambda(n)\ll x^{3/4}$. Hence, from (\ref{app}) and (\ref{app2}), we conclude that
	
	$$\sum_{\substack{n\leq x^{1/4}\\P^-(n)>z}}\!\!\!f_{\Gamma}(n)S(x/n,\Lambda_{(f_{\Gamma})_z})\ll_{\Gamma}\,\frac{x}{u}+\frac{x}{u\log x}\!\!\sum_{\substack{n\leq x^{3/4}\\P^-(n)>z}}\!\!\!\frac{|\Lambda_{f_{\Gamma}}(n)|}{n}.$$
	
	\vspace{1mm}
	
	\noindent
     The sum of the right-hand side is $\leq\sum_{n\leq x^{3/4}}|\Lambda_{f_{\Gamma}}(n)|/n\leq (D+m)\sum_{n\leq x^{3/4}}\Lambda(n)/n\ll \log x$, and so
	
	$$\sum_{\substack{n\leq x^{1/4}\\P^-(n)>z}}\!\!\!f_{\Gamma}(n)S(x/n,\Lambda_{(f_{\Gamma})_z})\ll_{\Gamma}\frac{x}{u}.\,\,\,$$
	
	\vspace{1mm}
	
	\noindent
	In addition, for $v\geq 1$ we have that
	
	\begin{eqnarray*}
		\,\,\,\,\,\,\,|S(v,\Lambda_{(f_{\Gamma})_z})-S(v,\Lambda_{f_{\Gamma}})|\leq\sum_{\substack{p\leq z\\p^{\nu}\leq v}}|\Lambda_{f_{\Gamma}}(p^{\nu})|\leq (D+m)\sum_{\substack{p\leq z\\p^{\nu}\leq v}}\Lambda(p^{\nu})\ll\log v\sum_{p\leq z}1\ll z\log v.
	\end{eqnarray*}

\vspace{1mm}
	
	\noindent
	Consequently, since $|f_{\Gamma}|\leq \tau_{D+m}$ and $\sum_{n\leq y}\tau_{D+m}(n)/n\ll (\log y)^{D+m}$ for any $y\geq 1$, it follows that
	
	\begin{eqnarray}\label{res}
	S(x,{\Lambda_{f_{\Gamma}}})\!\!\!&+&\!\!\!\!\!\sum_{\substack{1<n\leq x^{1/4}\\P^-(n)>z\,}}\!\!\!f_{\Gamma}(n)S(x/n,\Lambda_{f_{\Gamma}})=\!\!\!\sum_{\substack{n\leq x^{1/4}\\P^-(n)>z}}\!\!\!f_{\Gamma}(n)S(x/n,\Lambda_{f_{\Gamma}})\nonumber \\&=&\!\!\!\sum_{\substack{n\leq x^{1/4}\\P^-(n)>z}}\!\!\!f_{\Gamma}(n)S(x/n,\Lambda_{(f_{\Gamma})_z})+O(zx^{1/4}(\log x)^{D+m})\ll_{\Gamma} \frac{x}{u}.\,\,\,\,\,
	\end{eqnarray}
	
	\vspace{1mm}
	
	\noindent
	For the last estimate we used the fact that $z\leq x^{1/8}$, an inequality which follows from the choice of $z$ that we made. The Brun-Titchmarsh theorem guarantees that
	
	$$S(x/n,\Lambda_{f_{\Gamma}})-S(t/n,\Lambda_{f_{\Gamma}})\ll \frac{x}{nu}$$
	
	\vspace{2mm}
	
	\noindent
	for $t\in [x-x/u,x]$ and $n\leq x^{1/4}$. Using again the inequality $|f_\Gamma|\leq\tau_{D+m}$, we have
	
	\begin{eqnarray*}
		\sum_{\substack{n\leq x^{1/4}\\P^-(n)>z}}\!\!\!\frac{|f_{\Gamma}(n)|}{n}\ll(\log u)^{D+m}\!.                                    
	\end{eqnarray*}
	
	\noindent
	Thus, if we set $\Delta=x/u$, relation (\ref{res}) becomes
	
	\begin{eqnarray*}
		S(x,\Lambda_{f_{\Gamma}})\!\!\!&=&\!\!\!-\frac{1}{\Delta}\sum_{\substack{1<n\leq x^{1/4}\\P^-(n)>z}}\!\!\!f_{\Gamma}(n)\int_{x-\Delta}^xS(t/n,\Lambda_{f_{\Gamma}})\mathrm{d}t+O_{\Gamma}\!\left(\frac{x(\log u)^{D+m}}{u}\right)\\
		&=&\!\!\!-\frac{1}{\Delta}\sum_{\substack{1<n\leq x^{1/4}\\P^-(n)>z}}\!\!\!f_{\Gamma}(n)n\int_{\frac{x-\Delta}{n}}^{\frac{x}{n}}S(t,\Lambda_{f_{\Gamma}})\mathrm{d}t+O_{\Gamma}\!\left(\frac{x(\log u)^{D+m}}{u}\right)\\
		&=& \!\!\!-\frac{1}{\Delta}\int_{\frac{x-\Delta}{x^{1/4}}}^{\frac{x}{z}}S(t,\Lambda_{f_{\Gamma}})\Bigg(\sum_{\substack{(x-\Delta)/t<n\leq x/t\\P^-(n)>z}}\!\!\!\!\!\!f_{\Gamma}(n)n\Bigg)\mathrm{d}t+O_{\Gamma}\!\left(\frac{x(\log u)^{D+m}}{u}\right)\!\!.\,\,
	\end{eqnarray*}

\vspace{3mm}
	
	\noindent
	Using Shiu's theorem \cite[Theorem 20.3, p.~209]{dimb} with the non-negative, divisor-bounded, multiplicative function $\tau_{D+m}\cdot\mathds{1}_{P^-(\cdot)>z}$, we obtain
	
	$$\Bigg|\sum_{\substack{(x-\Delta)/t<n\leq x/t\\P^-(n)>z}}\!\!\!\!\!\!f_{\Gamma}(n)n\Bigg|\leq \frac{x}{t}\sum_{\substack{(x-\Delta)/t<n\leq x/t\\P^-(n)>z}}\!\!\!\!\!\!\tau_{D+m}(n)\ll \frac{x\Delta(\log u)^{D+m}}{t^2\log x},$$
	
	\vspace{2mm}
	
	\noindent
	for $t\leq x/z$ and $\Delta=x/u$. Thus, we arrive at the estimate
	
	\begin{equation}\label{f1}
	\,\,\,\,S(x,\Lambda_{f_{\Gamma}})\ll \frac{x(\log u)^{D+m}}{\log x}\!\!\int_1^x\frac{|S(t,\Lambda_{f_{\Gamma}})|}{t^2}\mathrm{d}t+O_{\Gamma}\!\left(\frac{x(\log u)^{D+m}}{u}\right)\!\!.
	\end{equation}
	
	\vspace{4mm}
	
	We continue by bounding the integral of the right-hand side and we start with the Cauchy-Schwarz inequality which implies that
	
	\begin{eqnarray}\label{rec}
	\int_1^x\frac{|S(t,\Lambda_{f_{\Gamma}})|}{t^2}\mathrm{d}t\!\!\!&\leq&\!\!\!\left(\log x\int_1^x\frac{|S(t,\Lambda_{f_{\Gamma}})|^2}{t^3}\mathrm{d}t\right)^{1/2}\nonumber \\&\ll&\!\!\!\left(\log x\int_1^{+\infty}\frac{|S(t,\Lambda_{f_{\Gamma}})|^2}{t^{3+2/\log x}}\mathrm{d}t\right)^{1/2}\!\!.
	\end{eqnarray}
	
	\vspace{3mm}
	
	\noindent
	Partial summation and a suitable change of variables give
	
	$$s\int_0^{+\infty}\!S(e^u,\Lambda_{f_{\Gamma}})e^{-u\sigma}e^{-iut}\mathrm{d}u=-\frac{L'}{L}(s,f_{\Gamma}) \quad \text{when} \quad \!\sigma>1.\,\,\,$$
	
	\vspace{2mm}
	
	\noindent
	So, for $c=1+1/\log x$, Parseval's theorem allows us to write
	
	$$\,\,\,\int_1^{+\infty}\frac{|S(t,\Lambda_{f_{\Gamma}})|^2}{t^{3+2/\log x}}\mathrm{d}t=\int_0^{+\infty}|S(e^u,\Lambda_{f_{\Gamma}})|^2e^{-2(1+1/\log x)u}\mathrm{d}u=\frac{1}{2\pi}\int_{\mathbb{R}}\left|\frac{L'}{L}\left(c+it,f_{\Gamma}\right)\right|^2\frac{\mathrm{d}t}{c^2+t^2}$$
	
	\vspace{3mm}
	
	\noindent
	and (\ref{rec}) implies that
	
	\begin{equation}\label{f2}
	\,\,\int_1^x\frac{|S(t,\Lambda_{f_{\Gamma}})|}{t^2}\mathrm{d}t\ll \sqrt{\log x}\left(\int_{\mathbb{R}}\left|\frac{L'}{L}\left(c+it,f_{\Gamma}\right)\right|^2\frac{\mathrm{d}t}{c^2+t^2}\right)^{1/2}\!\!.
	\end{equation}
	
	\vspace{3mm}
	
	\noindent
	As was explained in Section 2, we proceed to splitting the integral into two parts. In the first part we integrate over the interval $[-T,T]$, whereas the range of integration of the second integral consists of the large values $|t|>T$. In the beginning of the proof we assumed that the largest multiplicity of the zeros of $L(\cdot,f)$ is at most $D-1$. So, using Proposition \ref{dsp} (c) and the non-vanishing of $L(s,f_{\Gamma})$ on the vertical line $\sigma=1$, we conclude that 
	
	\begin{equation}\label{final1}
	\int_{|t|\leq T}\left|\frac{L'}{L}\left(c+it,f_{\Gamma}\right)\right|^2\frac{\mathrm{d}t}{c^2+t^2}\ll \max_{|t|\leq T,\,\sigma\in[1,2]}\left|\frac{L'}{L}\left(\sigma+it,f_{\Gamma}\right)\right|^2\!\!=:C_f(T).\,\,\,\,
	\end{equation}
	
	\vspace{3mm}
	
	\noindent
    The first part has been estimated and it remains to bound the second one. Since $|\Lambda_{f_{\Gamma}}(n)n^{-ik}|\leq (D+m)\Lambda(n)$ for any two positive integers $n$ and $k$ and $\zeta(s)\asymp (s-1)^{-1}$ in a fixed region around $1$, Lemma \ref{mvt} yields that
	
	\begin{eqnarray}\label{final2}
	\int_{|t|>T}\left|\frac{L'}{L}\left(c+it,f_{\Gamma}\right)\right|^2\frac{\mathrm{d}t}{c^2+t^2}\!\!\!&\leq&\!\!\!\sum_{|k|> T-1/2}\int_{k-1/2}^{k+1/2}\left|\frac{L'}{L}\left(c+it,f_{\Gamma}\right)\right|^2\frac{\mathrm{d}t}{c^2+t^2}\nonumber \\
	&\leq&\!\!\!4\!\sum_{|k|> T/2}\frac{1}{k^2}\int_{-\frac{1}{2}}^{\frac{1}{2}}\left|\frac{L'}{L}\left(c+i(t+k),f_{\Gamma}\right)\right|^2\!\mathrm{d}t\nonumber \\
	&\ll&\!\!\!\Bigg(\sum_{k>T/2}\frac{1}{k^2}\Bigg)\!\cdot \!\int_{-\frac{1}{2}}^{\frac{1}{2}}\left|\frac{\zeta'}{\zeta}(c+it)\right|^2\!\mathrm{d}t\nonumber \\
	&\asymp&\!\!\!\frac{1}{T}\int_{-\frac{1}{2}}^{\frac{1}{2}}\frac{\mathrm{d}t}{\frac{1}{(\log x)^2}+t^2}\leq\frac{\log x}{T}\int_{\R}\frac{\mathrm{d}\alpha}{1+{\alpha}^2}\ll \frac{\log x}{T}.
	\end{eqnarray}
	
	\vspace{2mm}
	
	\noindent
	Combining (\ref{final1}) and (\ref{final2}), we get that
	
	$$\int_{\R}\left|\frac{L'}{L}\left(c+it,f_{\Gamma}\right)\right|^2\frac{\mathrm{d}t}{c^2+t^2}\ll C_f(T)+ \frac{\log x}{T}.$$
	
	\vspace{3mm}
	
	\noindent
	We now insert this estimate into (\ref{f2}) to find that 
	
	$$\int_1^x\frac{|S(t,\Lambda_{f_{\Gamma}})|}{t^2}\mathrm{d}t\ll \sqrt{C_f(T)\log x}+\frac{\log x}{\sqrt{T}}.$$
	
	\vspace{2mm}
	
	\noindent
	Together with (\ref{f1}), this implies that
	
	\begin{equation*}
	S(x,\Lambda_{f_{\Gamma}})\leq O_{f,T}\left(\frac{x(\log u)^{D+m}}{\sqrt{\log x}}\right)\!+O_{\Gamma}\!\left(\frac{x(\log u)^{D+m}}{\min\{u,\sqrt{T}\}}\right)\!\!.
	\end{equation*}
	
	\vspace{3mm}
	
	\noindent
	Recalling that we chose $u=\min\{(\log x)^{A-D-1},\sqrt{T}\}$ and that $\sum_{p\leq x}(f(p)+\sum_{\gamma\in\Gamma}p^{i\gamma})\log p=\sum_{p\leq x}f_{\Gamma}(p)\log p=S(x,\Lambda_{f_{\Gamma}})+O(\sqrt{x})$, the proof of Theorem \ref{mthm} readily follows.
\vspace{2mm}


\bibliographystyle{alpha}

\end{document}